\documentclass{amsart}
\usepackage{amsmath, amssymb, latexsym, graphicx, epstopdf}
\usepackage{euscript}
\usepackage{subfigure}
\usepackage{epsfig}
\usepackage{psfrag}
\usepackage{amsthm}
\usepackage{setspace}
\theoremstyle{plain}

\newtheorem{thm}{Theorem}[section]
\newtheorem*{thm*}{Theorem}
\newtheorem{lem}[thm]{Lemma}
\newtheorem*{lem*}{Lemma}

\newtheorem*{prop*}{Proposition}
\newtheorem{cor}[thm]{Corollary}
\newtheorem*{cor*}{Corollary}

\newtheorem*{cla*}{Claim}

\newtheorem*{cond*}{Condition}

\newtheorem*{fac*}{Fact}

\newtheorem*{que*}{Question}

\newtheorem*{prob*}{Problem}

\newtheorem*{con*}{Conjecture}

\newtheorem*{rem*}{Remark}

\newtheorem*{rems*}{Remarks}

\newtheorem*{defn*}{Definition}

\newcommand{\C}{\ensuremath{\mathbb{C}}}

\newcommand{\Z}{\ensuremath{\mathbb{Z}}}

\newcommand{\Q}{\ensuremath{\mathbb{Q}}}

\newcommand{\slc}{\mathrm{SL}_2\C}

\newcommand{\DD}{\mathcal{D}}

\newcommand{\brac}[2]{\left\langle \frac{#1}{#2} \right\rangle}
\newcommand{\sbrac}[1]{\left\langle #1 \right\rangle}

\begin{document}
\title{Denominators and differences of boundary slopes for (1,1)-knots}
\author{Jason Callahan}
\begin{abstract}
We show that every nonzero integer occurs in the denominator of a boundary slope for infinitely many (1,1)-knots and that infinitely many (1,1)-knots have boundary slopes of arbitrarily small difference. Specifically, we prove that for any integers $m,n>1$ with $n$ odd the exterior of the Montesinos knot $K(-1/2, m/(2m\pm1),1/n)$ in $S^3$ contains an essential surface with boundary slope $r = 2(n-1)^2/n$ if $m$ is even and $2(n+1)^2/n$ if $m$ is odd. If $n \geq 4m + 1$, we prove that $K(-1/2, m/(2m+1),1/n)$ also has a boundary slope whose difference with $r$ is $(8m-2)/(n^2-4mn+n)$, which decreases to 0 as $n$ increases. All of these knots are (1,1)-knots.
\end{abstract}
\maketitle

\section{Introduction}
Embedded essential (i.e., incompressible and $\partial$-incompressible) surfaces have long been crucial to the study of 3-manifolds. The most general strategy for constructing such surfaces in a 3-manifold $M$, established in \cite{CullerShalen83}, associates essential surfaces to ideal points of irreducible curves in the character variety equal to the complex algebraic set of characters of $\slc$-representations of $\pi_1(M)$. This technique is especially effective in identifying the boundary slopes of $M$, i.e., the elements of $\pi_1(\partial M)$ represented by boundary components of essential surfaces.

If $M$ is a knot exterior in $S^3$, then $\partial M$ is a torus and boundary slopes are identified with elements of $\Q \cup \{\infty\}$ so that a longitude has slope 0 and a meridian slope $\infty$. Boundary slopes have been computed for some classes of knots: 0 and $pq$ are the boundary slopes for the $(p,q)$-torus knot, only even integers occur as boundary slopes for two-bridge knots (\cite{HatcherThurston85}), and all rational numbers occur as boundary slopes for Montesinos knots (\cite{HatcherOertel89}).

In this paper, we show that every nonzero integer occurs in the denominator of a boundary slope for infinitely many (1,1)-knots, a class of knots that includes torus knots, two-bridge knots, and some other Montesinos knots and has attracted recent attention; e.g., \cite{EudaveMunoz06} and \cite{EudaveMunozetal09} investigate closed and meridional essential surfaces in (1,1)-knot exteriors. We also show that infinitely many (1,1)-knots have boundary slopes of arbitrarily small difference. Specifically, we review the algorithm of \cite{HatcherOertel89} in Section~\ref{prelim} and use it in Section~\ref{main} to prove:

\begin{thm}\label{res}
For any integers $m,n>1$ with $n$ odd, the exterior of the Montesinos knot $K(-1/2, m/(2m\pm1),1/n)$ in $S^3$ contains an essential surface with $n$ sheets, one boundary component, Euler characteristic $-n$, and boundary slope $2(n-1)^2/n$ if $m$ is even and $2(n+1)^2/n$ if $m$ is odd.
\end{thm}

\begin{thm}\label{res2}
For any integer $m\geq1$ and odd integer $n \geq 4m + 1$, the exterior of the Montesinos knot $K(-1/2, m/(2m+1),1/n)$ in $S^3$ contains an essential surface with $n-4m+1$ sheets, two boundary components, Euler characteristic $2m+3-n$, and boundary slope $2(n^2-4mn-n+8m-1)/(n-4m+1)$ if $m$ is even and $2(n^2-4mn+3n-8m+3)/(n-4m+1)$ if $m$ is odd.
\end{thm}

Thus, if $m>1$ and $n \geq 4m + 1$ as above, $K(-1/2, m/(2m+1),1/n)$ has boundary slopes with difference $(8m-2)/(n^2-4mn+n)$, which becomes arbitrarily small for any $m$ as $n$ increases. This generalizes the main result of \cite{Ichihara14}, which treats only the case $m=2$ but overlooks that the boundary slope is reducible by 2 and incorrectly concludes that the essential surface has only one boundary component and hence must be non-orientable; we explain this in Section~\ref{main}. A review of definitions and related results in Section~\ref{prelim} reveals that all knots above are (1,1)-knots, and we conclude in Section~\ref{disc} with discussion of a bound on the denominators of boundary slopes for Montesinos knots from \cite{IchiharaMizushima07}.

\section{Preliminaries}\label{prelim}
As in \cite{HatcherOertel89}, we consider Montesinos knots $K(p_1/q_1, \dots, p_n/q_n)$ obtained by connecting $n$ rational tangles of irreducible non-integral slopes $p_1/q_1, \dots, p_n/q_n$ in a simple cyclic pattern that yields a knot if just one $q_i$ is even or if all $q_i$ are odd and the number of odd $p_i$ is odd. Two-bridge knots are the Montesinos knots with $n<3$, and $(q_1, \dots, q_n)$-pretzel knots are the Montesinos knots $K(1/q_1, \dots, 1/q_n)$.

Following \cite{MorimotoSakumaYokota96}, the \textit{tunnel number} $t(K)$ of a knot $K$ in $S^3$ is the minimum number of mutually disjoint arcs $\{\tau_i\}$ properly embedded in the exterior of $K$ such that the exterior of $K \cup (\cup \tau_i)$ is a handlebody, and $K$ has a \textit{(g,b)-decomposition} if there is a genus $g$ Heegaard splitting $\{W_1, W_2\}$ of $S^3$ such that $K$ intersects each $W_i$ in a $b$-string trivial arc system; then $t(K) \leq g+b-1$, so $K$ has tunnel number one if it admits a (1,1)-decomposition, i.e., is a (1,1)-knot. The following reveals that the Montesinos knots $K(-1/2, m/(2m\pm1),1/n)$ with $n$ odd are (1,1)-knots (\cite[Theorem 2.2 and a closing remark]{MorimotoSakumaYokota96}). 

\begin{thm}\label{MSY}
A Montesinos knot $K(p_1/q_1, \dots, p_n/q_n)$ is a (1,1)-knot if and only if $n=2$ or $n=3$ and, up to cyclic permutation of the indices, either $q_1 = 2$ and $q_2$  and $q_3$ are odd or $p_2/q_2 \equiv p_3/q_3 \equiv \pm1/3 \emph{ mod }1$ and $\sum p_i/q_i = \pm1/(3q_1)$.
\end{thm}

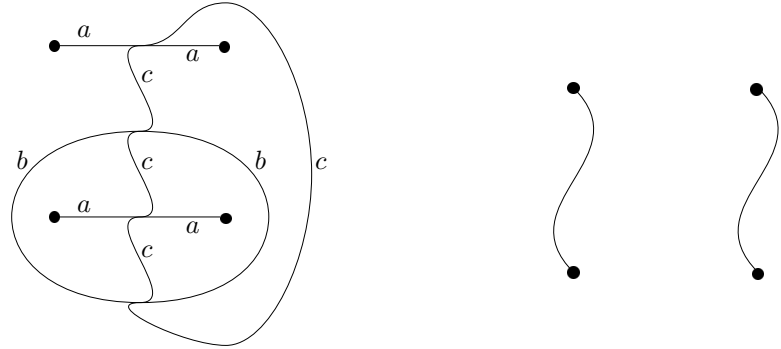
\begin{figure}[ht]
\begin{center}
\subfigure[The train track with projective weights $a$, $b$, and $c$.]{
\begin{picture}(200,125)
\put(25,0){\includegraphics[width=4cm]{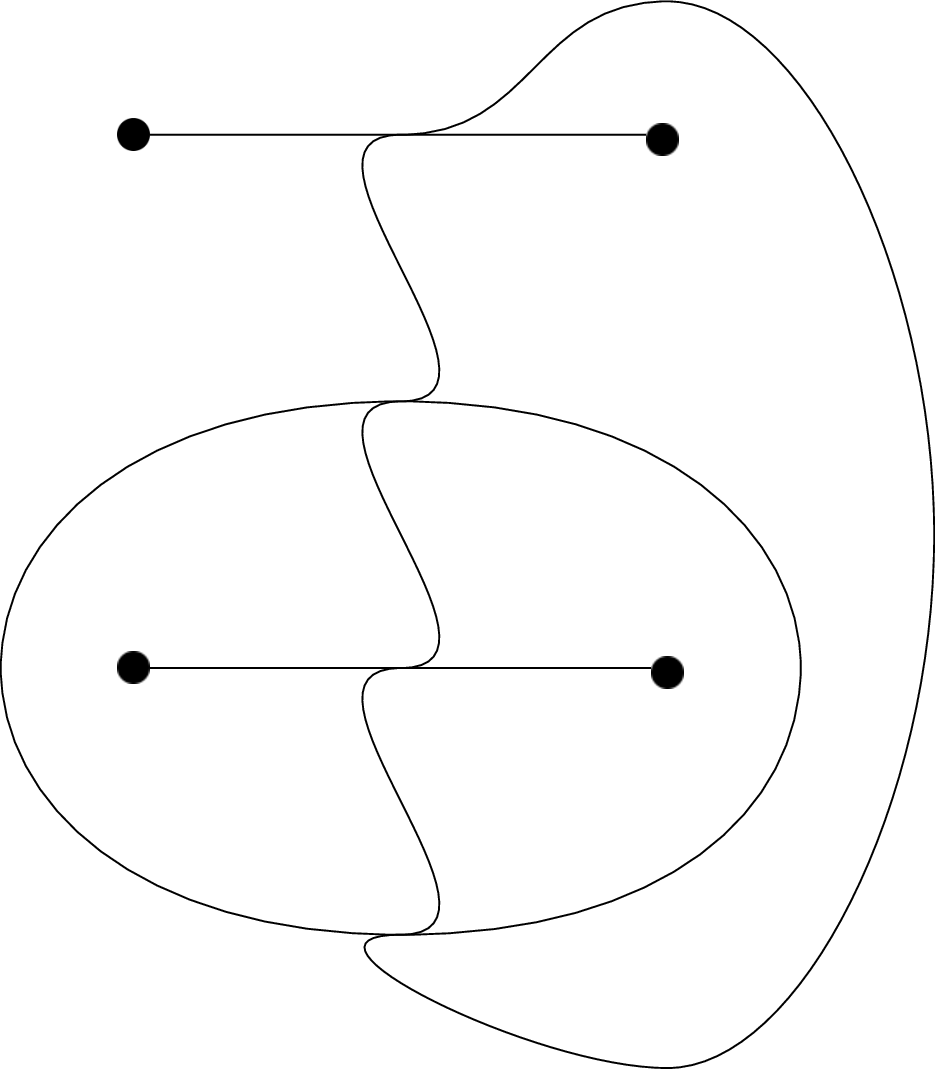}}
\put(50,117){$a$}
\put(91,108){$a$}
\put(91,43){$a$}
\put(50,51){$a$}
\put(74,100){$c$}
\put(74,67){$c$}
\put(140,67){$c$}
\put(74,34){$c$}
\put(27,67){$b$}
\put(117,67){$b$}
\end{picture}
\label{traintrack}}
\subfigure[The $\infty$-tangle.]{
\begin{picture}(100,100)
\put(20,25){\includegraphics[width=3cm]{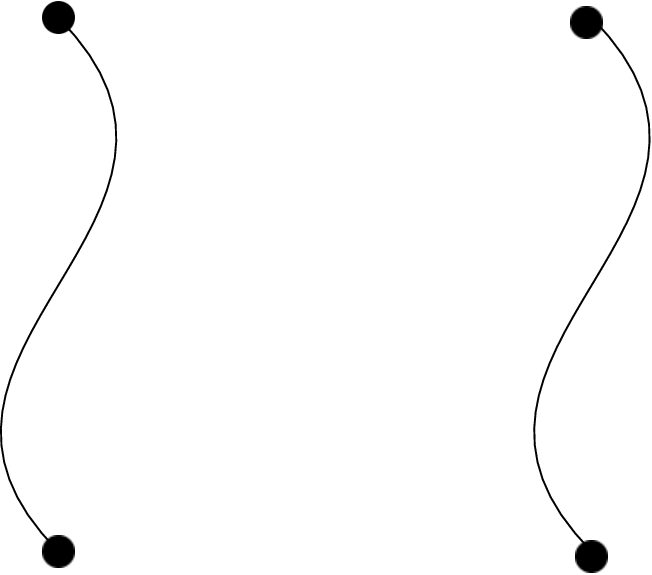}}
\end{picture}
\label{infty}}
\end{center}
\caption{Curve systems and tangles.}  
\end{figure}

We now recall the algorithm of \cite{HatcherOertel89} to compute boundary slopes for Montesinos knots $K=K(p_1/q_1, \dots p_n/q_n)$ with $n \geq 3$ as also presented in \cite{ChesebroTillman07}, \cite{IchiharaMizushima07}, \cite{Wu11}, and \cite{Ichihara14}. See those for details, but briefly \cite{HatcherOertel89} associates \textit{candidate surfaces} to \textit{admissible edgepath systems} in a graph $\DD$ in the $uv$-plane whose vertices $(u,v)$ correspond to projective curve systems $[a,b,c]$ on the 4-punctured sphere carried by the train track in Fig.~\ref{traintrack} via $u=b/(a+b)$ and $v=c/(a+b)$. Specifically, the vertices of $\DD$ are:
\begin{itemize}
\item the $\infty$-\textit{tangle} $\sbrac{\infty}$ in Fig.~\ref{infty} with $uv$-coordinates $(-1,0)$, 
\item the $p/q$-\textit{circles} $\sbrac{p/q}^\circ$ whose $uv$-coordinates $(1,p/q)$ correspond to the projective curve systems $[0,q,p]$, and
\item the $p/q$-\textit{tangles} $\sbrac{p/q}$ whose $uv$-coordinates $((q-1)/q,p/q)$ correspond to the projective curve systems $[1,q-1,p]$.
\end{itemize}

\begin{figure}[ht] 
\begin{center}
\begin{picture}(175, 215)
  \put(0,0){\includegraphics[width=6cm]{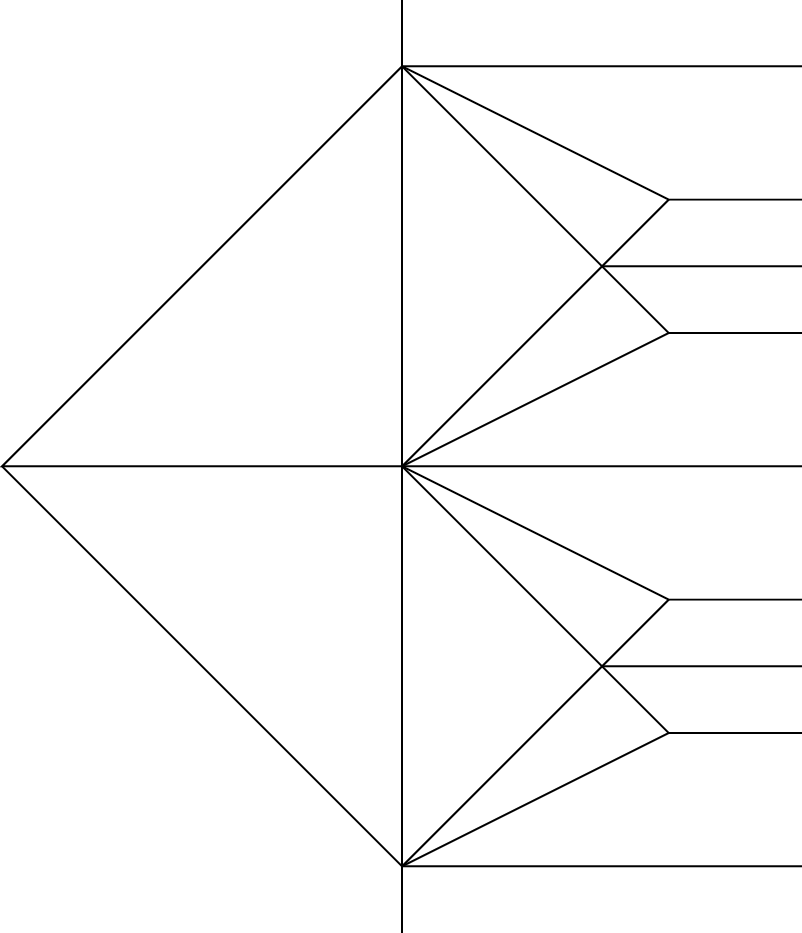}}
  \put(-18,97){$\sbrac{\infty}$}
  \put(71,185){$\sbrac{1}$}
  \put(171,182){$\sbrac{1}^\circ$}
 \put(171,154){$\sbrac{2/3}^\circ$}
 \put(139,161){$\sbrac{2/3}$}
 \put(139,118){$\sbrac{1/3}$}
 \put(171,140){$\sbrac{1/2}^\circ$}
 \put(101,140){$\sbrac{1/2}$}
 \put(171,126){$\sbrac{1/3}^\circ$}
  \put(71,104){$\sbrac{0}$}
\put(171,97){$\sbrac{0}^\circ$}
 \put(171,40){$\sbrac{-2/3}^\circ$}
 \put(139,33){$\sbrac{-2/3}$}
 \put(139,76){$\sbrac{-1/3}$}
 \put(171,55){$\sbrac{-1/2}^\circ$}
 \put(93,54){$\sbrac{-1/2}$}
 \put(171,69){$\sbrac{-1/3}^\circ$}
  \put(64,9){$\sbrac{-1}$}
 \put(171,12){$\sbrac{-1}^\circ$}
  \put(84,202){$\vdots$}
  \put(84,-11){$\vdots$}
  \end{picture}
\end{center}
\caption{Part of the graph $\DD$.} \label{graph}
\end{figure}

If $|ps-qr|=1$, then $\left[\sbrac{p/q}, \sbrac{r/s}\right]$ is a \textit{non-horizontal edge} in $\DD$ connecting $\sbrac{r/s}$ to $\sbrac{p/q}$; the remaining edges in $\DD$ are:
\begin{itemize}
\item the \textit{horizontal edges} $\left[\sbrac{p/q}, \sbrac{p/q}^\circ\right]$ connecting $\sbrac{p/q}^\circ$ to $\sbrac{p/q}$,
\item the \textit{vertical edges} $\left[\sbrac{m}, \sbrac{m+1}\right]$ connecting $\sbrac{m+1}$ to $\sbrac{m}$, and
\item the \textit{infinity edges} $\left[\sbrac{\infty}, \sbrac{m}\right]$ connecting $\sbrac{m}$ to $\sbrac{\infty}$
\end{itemize}
for any integer $m$. Fig.~\ref{graph} shows part of the graph $\DD$.

\textit{Rational points} $(p/q, r/s) \in \DD \cap \Q^2$ need not be vertices of $\DD$ and correspond to the projective curve systems $[s(q-p),sp,rq]$. If $\left[\sbrac{p/q}, \sbrac{r/s}\right]$ is a non-horizontal edge in $\DD$, then $\frac{k}{m}\sbrac{p/q} + \frac{m-k}{m}\sbrac{r/s}$ is a rational point on this edge with coordinates
\begin{equation}\label{coord}
\left(\frac{k(q-s)+m(s-1)}{k(q-s)+ms},\frac{k(p-r)+mr}{k(q-s)+ms}\right)
\end{equation}

If $\left[\sbrac{p/q}, \sbrac{p/q}^\circ\right]$ is a horizontal edge in $\DD$, then $\frac{k}{m}\sbrac{p/q} + \frac{m-k}{m}\sbrac{p/q}^\circ$ is a rational point on this edge with coordinates $((mq-k)/mq,p/q)$.

An \textit{edgepath} in $\DD$ is a piecewise linear path $[0,1] \to \DD$ that begins and ends at rational points (not necessarily vertices) of $\DD$. An \textit{admissible edgepath system} $\gamma=(\gamma_1, \ldots, \gamma_n)$ is an $n$-tuple of edgepaths in $\DD$ such that:
\begin{itemize}
\item[(E1)] Each starting point $\gamma_i(0)$ lies on the horizontal edge $\left[\sbrac{p_i/q_i}, \sbrac{p_i/q_i}^\circ\right]$, and $\gamma_i$ is constant if $\gamma_i(0) \neq \sbrac{p_i/q_i}$.
\item[(E2)] Each $\gamma_i$ is \textit{minimal}, i.e., it never stops and retraces itself, and it never travels along two sides of a triangle in $\DD$ in succession.
\item[(E3)] The ending points $\gamma_1(1), \dots, \gamma_n(1)$ all lie on a vertical line (i.e., have the same $u$-coordinate), and their $v$-coordinates sum to zero.
\item[(E4)] Each $\gamma_i$ proceeds monotonically from right to left where traversing vertical edges is permitted, i.e., if $0\leq t_1 < t_2 \leq 1$, then the $u$-coordinate of $\gamma_i(t_1)$ is at least as great as the $u$-coordinate of $\gamma_i(t_2)$.
\end{itemize}  

The aforementioned curve systems on the 4-punctured sphere describe how the boundaries of 3-balls that decompose $S^3$ and each contain a tangle of $K$ intersect properly embedded surfaces in the exterior of $K$, and \cite{HatcherOertel89} associates a finite number of \textit{candidate surfaces} to each admissible edgepath system; every essential surface in the exterior of $K$ with non-empty boundary of finite slope is then isotopic to one of the candidate surfaces (\cite[Proposition 1.1]{HatcherOertel89}). Conversely, we have the following (\cite[Lemma 3]{Ichihara14} and \cite[Proposition 2.1]{HatcherOertel89}).

\begin{lem}\label{ess}
If an admissible edgepath system has ending points with positive $u$-coordinate and final edges all traveling in the same direction (upward or downward) from right to left, then all candidate surfaces are essential.
\end{lem}

\begin{lem}\label{const}
A candidate surface is essential if any of its edgepaths is constant.
\end{lem}

A candidate surface $S$ for an admissible edgepath system $\gamma = (\gamma_1, \dots, \gamma_n)$ has $\sharp S$ \textit{sheets} if it meets a small meridional circle of $K$ minimally in $\sharp S$ points. To compute this as in \cite{IchiharaMizushima07} and \cite{Wu11}, define the length $|\gamma_i|$ of each $\gamma_i$ by counting the length of a full edge as 1 and the length of a partial edge from $\sbrac{r/s}$ to $\frac{k}{m}\sbrac{p/q} + \frac{m-k}{m}\sbrac{r/s}$ as $k/m$; constant edges have length 0, and the length of $\gamma$ is $|\gamma| = \sum|\gamma_i|$. If $\gamma_i$ is not constant, let $l_i$ be the least positive integer so that $l_i|\gamma_i| \in \Z$. If $\gamma_i$ is the constant edgepath $[\frac{k}{m}\sbrac{p_i/q_i} + \frac{m-k}{m}\sbrac{p_i/q_i}^\circ]$, let $l_i$ be the least positive integer so that $l_im/k \in \Z$. Then $\sharp S  = \text{lcm}(l_1, \dots, l_n)$, and the Euler characteristic of $S$ is (\cite[Formula 3.5]{IchiharaMizushima07}):
\begin{equation}\label{euler}
\chi(S) = -\sharp S\left(|\gamma| + n_c - n + \left(n - 2 - \sum_{\text{const}} \frac{1}{q_i}\right)\frac{1}{1-u}\right)
\end{equation}
where $n_c$ is the number of constant edgepaths, the sum is over constant edgepaths $\gamma_i = [\frac{k}{m}\sbrac{p_i/q_i} + \frac{m-k}{m}\sbrac{p_i/q_i}^\circ]$, and $u$ is the $u$-coordinate of the ending points $\gamma_1(1), \dots, \gamma_n(1)$.

To compute the boundary slope of a candidate surface $S$ for an admissible edgepath system $\gamma$ as in \cite{HatcherOertel89}, define the \textit{twist number} of $S$ as $\tau(S) = 2(e_{-}-e_{+})$, where $e_{+}$ ($e_{-}$) is the number of edges of $\gamma$ that travel upward (downward) from right to left (infinity edges are not counted). Fractional values of $e_{\pm}$ correspond to edges of $\gamma$ that only traverse a fraction of an edge in $\DD$; the segment from $\sbrac{r/s}$ to $\frac{k}{m}\sbrac{p/q} + \frac{m-k}{m}\sbrac{r/s}$ counts as $k/m$ of an edge. The boundary slope of $S$ is $\tau(S)-\tau(\Sigma)$, where $\Sigma$ is a Seifert surface for $K$ that is a candidate surface found as follows (\cite[pages 460-461]{HatcherOertel89}).

\begin{lem}\label{seif}
A candidate surface associated to an admissible edgepath system $\gamma=(\gamma_1, \ldots, \gamma_n)$ is a Seifert surface for $K(p_1/q_1, \dots, p_n/q_n)$ if one $q_i$ is even and each $\gamma_i$ is a minimal edgepath from $\sbrac{p_i/q_i}$ to $\sbrac{\infty}$ whose\emph{ mod 2 }reduction uses only one edge of the triangle in $\DD$ with vertices $\sbrac{\infty}$, $\sbrac{0}$, and $\sbrac{1}$ such that the total number of odd-integer vertices in $\gamma$ is even.
\end{lem}

The denominator of a boundary slope $p/q$ in lowest terms is the minimum number of points each boundary component of an essential surface $S$ with $\sharp S$ sheets meets a small meridional circle of $K$, so the number of boundary components of $S$ is $|\partial S| = \sharp S/q$.

Finally, as noted in \cite{HatcherOertel89}, a candidate surface may be non-orientable and incompressible only in the weaker sense, not $\pi_1$-injective (\cite{HatcherThurston85} also uses this weaker definition of incompressibility), but an orientable candidate surface associated to the same edgepath system (and hence same boundary slope) also exists (but possibly with twice as many sheets) as follows (\cite[Lemma 2.2]{Wu11}).

\begin{lem}\label{orient}
If $S$ is a candidate surface with $\sharp S$ sheets associated to an admissible edgepath system, then there exists an orientable candidate surface with $\sharp S$ or  $2\sharp S$ sheets associated to the same edgepath system.
\end{lem}

\section{Our Results}\label{main}
We now prove Theorem \ref{res} restated here for convenience.

\begin{thm}\label{reres}
For any integers $m,n>1$ with $n$ odd, the exterior of the Montesinos knot $K=K(-1/2, m/(2m\pm1),1/n)$ contains an essential surface with $n$ sheets, one boundary component, Euler characteristic $-n$, and boundary slope $2(n-1)^2/n$ if $m$ is even and $2(n+1)^2/n$ if $m$ is odd.
\end{thm}

\begin{proof}
Let $\gamma$ be the edgepath system given by
\begin{eqnarray*}
\gamma_1 &=& \left[\frac{1}{n}\sbrac{-1} + \frac{n-1}{n}\sbrac{-\frac{1}{2}}, \sbrac{-\frac{1}{2}}\right]\\
\gamma_2 &=& \left[\frac{1}{n}\sbrac{0} + \frac{n-1}{n}\sbrac{\frac{1}{2}}, \sbrac{\frac{1}{2}}, \sbrac{\frac{m}{2m\pm1}}\right]\\
\gamma_3 &=& \left[\frac{n-1}{n}\sbrac{0} + \frac{1}{n}\brac{1}{n}, \sbrac{\frac{1}{n}}\right]
\end{eqnarray*}

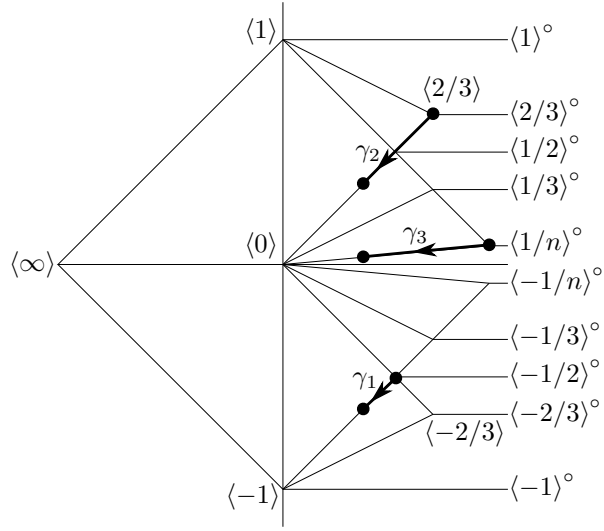
\begin{figure}[ht] 
\begin{center}
\begin{picture}(175, 200)
  \put(0,0){\includegraphics[width=6cm]{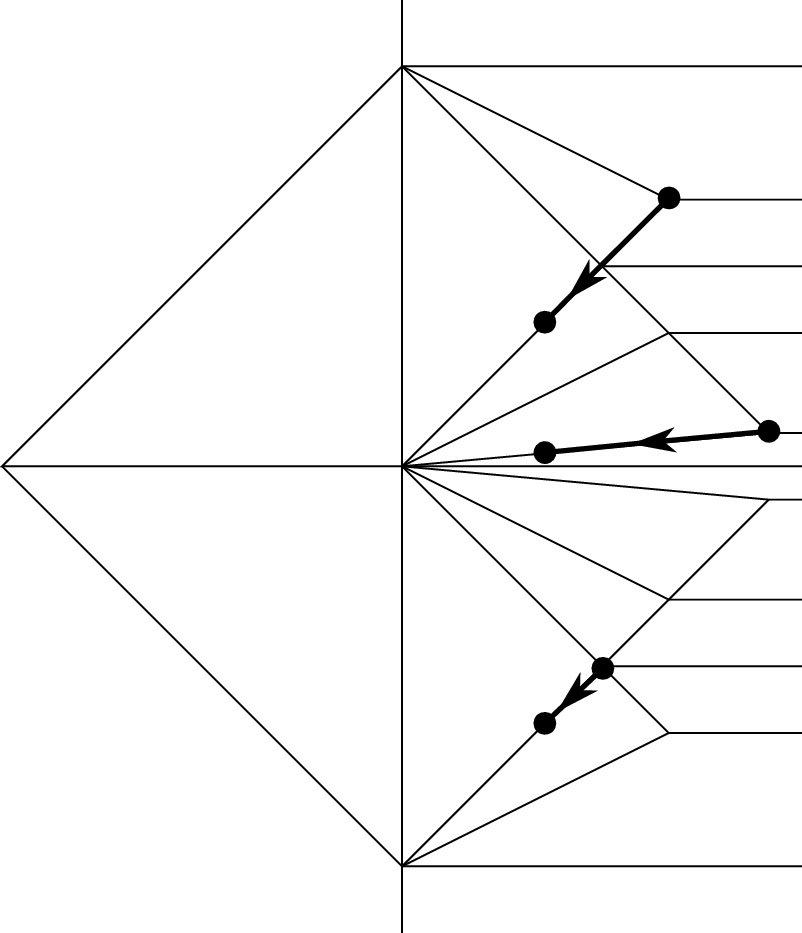}}
  \put(-18,97){$\sbrac{\infty}$}
  \put(71,185){$\sbrac{1}$}
  \put(171,182){$\sbrac{1}^\circ$}
 \put(171,154){$\sbrac{2/3}^\circ$}
 \put(138,163){$\sbrac{2/3}$}
 \put(171,140){$\sbrac{1/2}^\circ$}
 \put(113,140){$\gamma_2$}
 \put(112,54){$\gamma_1$}
 \put(131,109){$\gamma_3$}
 \put(171,126){$\sbrac{1/3}^\circ$}
  \put(71,104){$\sbrac{0}$}
\put(171,105){$\sbrac{1/n}^\circ$}
\put(171,90){$\sbrac{-1/n}^\circ$}
 \put(171,40){$\sbrac{-2/3}^\circ$}
 \put(139,33){$\sbrac{-2/3}$}
 \put(171,55){$\sbrac{-1/2}^\circ$}
 \put(171,69){$\sbrac{-1/3}^\circ$}
  \put(64,9){$\sbrac{-1}$}
 \put(171,12){$\sbrac{-1}^\circ$}
    \end{picture}
\end{center}
\caption{The edgepath system $\gamma$ for the case of $2m-1$ with $m=2$ and $n=3$.} \label{gamma}
\end{figure}

See Fig.~\ref{gamma} for the case of $2m-1$ with $m=2$ and $n=3$. To obtain an associated candidate surface, we verify that $\gamma$ satisfies conditions (E1-4):
\begin{itemize}
\item[(E1)] $\gamma_1(0) = \sbrac{-1/2}$ lies on $\left[\sbrac{-1/2}, \sbrac{-1/2}^\circ\right]$, $\gamma_2(0) = \sbrac{m/(2m\pm1)}$ lies on $\left[\sbrac{m/(2m\pm1)}, \sbrac{m/(2m\pm1)}^\circ\right]$, and $\gamma_3(0) = \sbrac{1/n}$ lies on $\left[\sbrac{1/n}, \sbrac{1/n}^\circ\right]$; none of the $\gamma_i$ are constant.
\item[(E2)] No $\gamma_i$ stops and retraces itself or travels along two sides of a triangle in $\DD$ in succession.
\item[(E3)] Using Formula (\ref{coord}),
\begin{eqnarray*}
\gamma_1(1) &=& \frac{1}{n}\sbrac{-1} + \frac{n-1}{n}\sbrac{-\frac{1}{2}} = \left(\frac{n-1}{2n-1},\frac{-n}{2n-1}\right)\\
\gamma_2(1) &=& \frac{1}{n}\sbrac{0} + \frac{n-1}{n}\sbrac{\frac{1}{2}} = \left(\frac{n-1}{2n-1},\frac{n-1}{2n-1}\right)\\
\gamma_3(1) &=& \frac{n-1}{n}\sbrac{0} + \frac{1}{n}\brac{1}{n} = \left(\frac{n-1}{2n-1},\frac{1}{2n-1}\right)
\end{eqnarray*}
These lie on a vertical line (i.e., have the same $u$-coordinate), and their $v$-coordinates sum to zero.
\item[(E4)] Each $\gamma_i$ proceeds monotonically from right to left.
\end{itemize}  

Hence, $\gamma$ is an admissible edgepath system, so a candidate surface $S$ in the exterior of $K$ can be associated to it. Because $n>1$, the $u$-coordinate of the ending points $\gamma_i(1)$ is positive, and all final edges travel downward from right to left; therefore, $S$ is essential by Lemma~\ref{ess}.

Now $|\gamma_1| = 1/n$, $|\gamma_2| = (n+1)/n$, and $|\gamma_3| = (n-1)/n$, so $n$ is the least positive integer such that $n|\gamma_i| \in \Z$ for all $i$, so $S$ has $n$ sheets and Euler characteristic $-n$ by Formula (\ref{euler}).

In the case of $2m-1$, all edges of $\gamma$ travel downward from right to left, so $e_+=0$, $e_-=(2n+1)/n$, and $\tau(S) = (4n+2)/n$.

In the case of $2m+1$, $\gamma_2$ is one edge upward followed by $1/n$ of an edge downward from right to left; all other edges remain downward from right to left, so $e_+=1$, $e_-=(n+1)/n$, and $\tau(S) = 2/n$.

We now find a Seifert surface for $K$ that is a candidate surface for an admissible edgepath system. If $m$ is even, let $\delta$ be the edgepath system
\begin{eqnarray*}
\delta_1 & = & \left[\sbrac{\infty}, \sbrac{-1}, \sbrac{-\frac{1}{2}}\right]\\
\delta_2 & = & \left[\sbrac{\infty}, \sbrac{0}, \brac{1}{2}, \brac{m}{2m\pm1}\right]\\
\delta_3 & = & \left[\sbrac{\infty}, \sbrac{1}, \brac{1}{2}, \dots, \brac{1}{n-1}, \brac{1}{n}\right]
\end{eqnarray*}

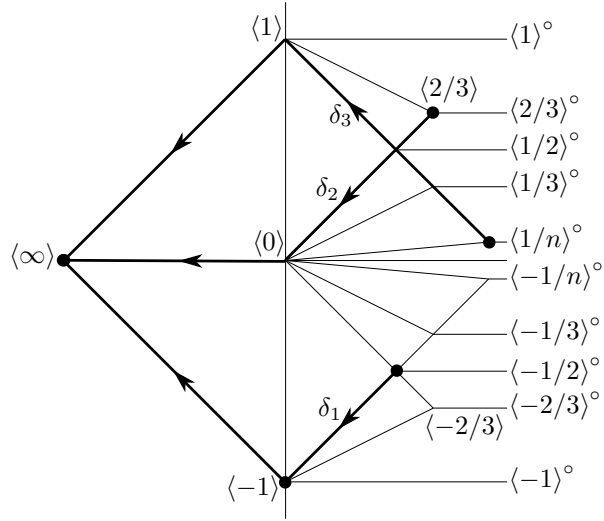
\begin{figure}[ht] 
\begin{center}
\begin{picture}(175, 200)
  \put(0,0){\includegraphics[width=6cm]{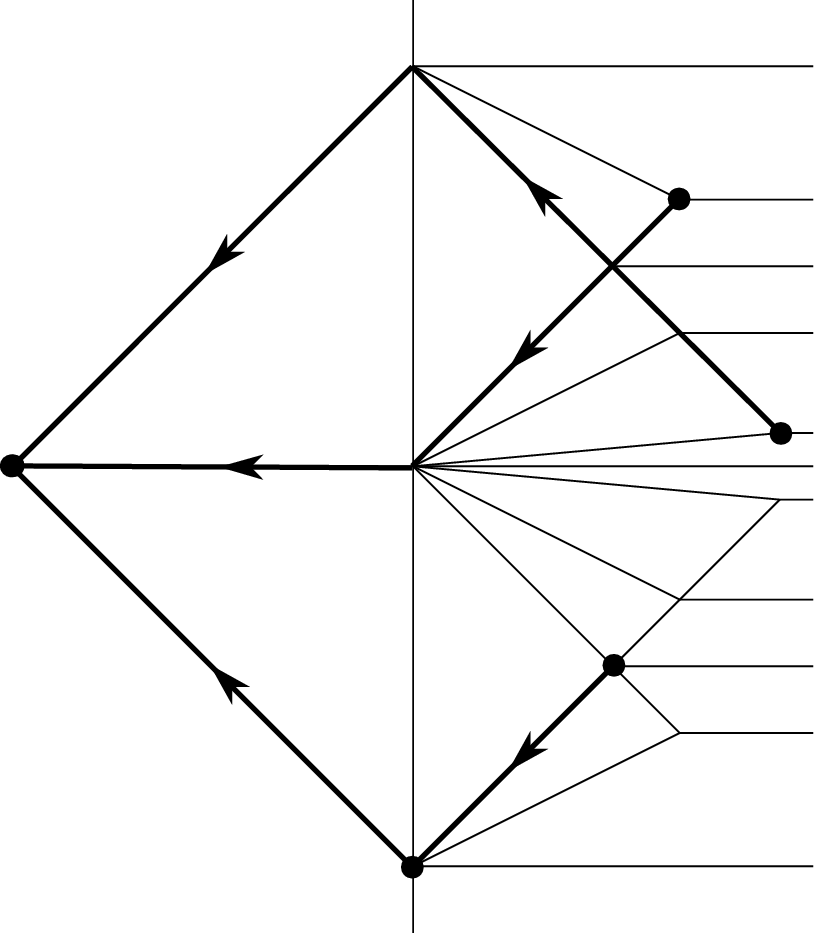}}
  \put(-18,97){$\sbrac{\infty}$}
  \put(73,183){$\sbrac{1}$}
  \put(171,180){$\sbrac{1}^\circ$}
 \put(171,152){$\sbrac{2/3}^\circ$}
 \put(137,160){$\sbrac{2/3}$}
 \put(171,138){$\sbrac{1/2}^\circ$}
 \put(98,122){$\delta_2$}
 \put(99,39){$\delta_1$}
 \put(103,150){$\delta_3$}
 \put(171,124){$\sbrac{1/3}^\circ$}
  \put(74,102){$\sbrac{0}$}
\put(171,103){$\sbrac{1/n}^\circ$}
\put(171,89){$\sbrac{-1/n}^\circ$}
 \put(171,40){$\sbrac{-2/3}^\circ$}
 \put(139,32){$\sbrac{-2/3}$}
 \put(171,53){$\sbrac{-1/2}^\circ$}
 \put(171,68){$\sbrac{-1/3}^\circ$}
  \put(64,9){$\sbrac{-1}$}
 \put(171,12){$\sbrac{-1}^\circ$}
  \end{picture}
\end{center}
\caption{The edgepath system $\delta$ for the case of $2m-1$ with $m=2$ and $n=3$.} \label{delta}
\end{figure}

See Fig.~\ref{delta} for the case of $2m-1$ with $m=2$ and $n=3$. To obtain an associated candidate surface, we verify that $\delta$ satisfies conditions (E1-4):
\begin{itemize}
\item[(E1)] $\delta_1(0) = \sbrac{-1/2}$ lies on $\left[\sbrac{-1/2}, \sbrac{-1/2}^\circ\right]$, $\delta_2(0) = \sbrac{m/(2m\pm1)}$ lies on $\left[\sbrac{m/(2m\pm1)}, \sbrac{m/(2m\pm1)}^\circ\right]$, and $\delta_3(0) = \sbrac{1/n}$ lies on $\left[\sbrac{1/n}, \sbrac{1/n}^\circ\right]$; none of the $\delta_i$ are constant.
\item[(E2)] No $\delta_i$ stops and retraces itself or travels along two sides of a triangle in $\DD$ in succession.
\item[(E3)] Each $\delta_i(1) = \sbrac{\infty} = (-1,0)$, so they all lie on a vertical line (i.e., have the same $u$-coordinate), and their $v$-coordinates sum to zero.
\item[(E4)] Each $\delta_i$ proceeds monotonically from right to left.
\end{itemize}  

Hence, $\delta$ is an admissible edgepath system with one $q_i$ even and each $\delta_i$ a minimal edgepath from $\sbrac{p_i/q_i}$ to $\sbrac{\infty}$ with integer vertices $\sbrac{-1}$, $\sbrac{0}$, and $\sbrac{1}$ and mod 2 reductions that use only the edges $[\sbrac{\infty}, \sbrac{1}]$, $[\sbrac{\infty}, \sbrac{0}]$, and $[\sbrac{\infty}, \sbrac{1}]$ respectively, so a candidate surface $\Sigma$ is a Seifert surface for $K$ by Lemma \ref{seif}.

Ignoring the infinity edges, in the case of $2m-1$, $\delta_3$ consists of $n-1$ edges traveling upward from right to left; all other edges travel downward from right to left, so $e_+=n-1$, $e_-=3$, $\tau(\Sigma) = 8-2n$, and the boundary slope of $S$ is $\tau(S) - \tau(\Sigma) = 2(n-1)^2/n$.

Again ignoring the infinity edges, in the case of $2m+1$, $\delta_2$ is one edge traveling upward followed by one traveling downward from right to left; all other edges remain unchanged, so $e_+=n$, $e_-=2$, $\tau(\Sigma) = 4-2n$, and the boundary slope of $S$ is again $\tau(S) - \tau(\Sigma) = 2(n-1)^2/n$.

If $m$ is odd, let $\delta$ be the edgepath system
\begin{eqnarray*}
\delta_1 & = & \left[\sbrac{\infty}, \sbrac{0}, \sbrac{-\frac{1}{2}}\right]\\
\delta_2 & = & \left[\sbrac{\infty}, \sbrac{1}, \brac{1}{2}, \brac{m}{2m\pm1}\right]\\
\delta_3 & = & \left[\sbrac{\infty}, \sbrac{1}, \brac{1}{2}, \dots, \brac{1}{n-1}, \brac{1}{n}\right]
\end{eqnarray*}

We again verify that $\delta$ satisfies conditions (E1-4):
\begin{itemize}
\item[(E1)] $\delta_1(0) = \sbrac{-1/2}$ lies on $\left[\sbrac{-1/2}, \sbrac{-1/2}^\circ\right]$, $\delta_2(0) = \sbrac{m/(2m\pm1)}$ lies on $\left[\sbrac{m/(2m\pm1)}, \sbrac{m/(2m\pm1)}^\circ\right]$, and $\delta_3(0) = \sbrac{1/n}$ lies on $\left[\sbrac{1/n}, \sbrac{1/n}^\circ\right]$; none of the $\delta_i$ are constant.
\item[(E2)] No $\delta_i$ stops and retraces itself or travels along two sides of a triangle in $\DD$ in succession.
\item[(E3)] Each $\delta_i(1) = \sbrac{\infty} = (-1,0)$, so they all lie on a vertical line (i.e., have the same $u$-coordinates), and their $v$-coordinates sum to zero.
\item[(E4)] Each $\delta_i$ proceeds monotonically from right to left.
\end{itemize}  

Hence, $\delta$ is an admissible edgepath system with one $q_i$ even and each $\delta_i$ a minimal edgepath from $\sbrac{p_i/q_i}$ to $\sbrac{\infty}$ with integer vertices $\sbrac{0}$, $\sbrac{1}$, and $\sbrac{1}$ and mod 2 reductions that use only the edges $[\sbrac{\infty}, \sbrac{0}]$, $[\sbrac{\infty}, \sbrac{1}]$, and $[\sbrac{\infty}, \sbrac{1}]$ respectively, so a candidate surface $\Sigma$ is a Seifert surface for $K$ by Lemma \ref{seif}.

Ignoring the infinity edges, in the case of $2m-1$, $\delta_2$ is one edge traveling downward followed by one traveling upward from right to left; all other edges travel upward from right to left, so $e_+=n+1$, $e_-=1$, $\tau(\Sigma) = -2n$, and the boundary slope of $S$ is $\tau(S) - \tau(\Sigma) = 2(n+1)^2/n$.

Again ignoring the infinity edges, in the case of $2m+1$, all edges travel upward from right to left, so $e_+=n+2$, $e_-=0$, $\tau(\Sigma) = -4-2n$, and the boundary slope of $S$ is again $\tau(S) - \tau(\Sigma) = 2(n+1)^2/n$.

Because $n$ is odd, the boundary slope $2(n\pm1)^2/n$ is in lowest terms, so since $S$ has $n$ sheets, it has one boundary component for $m$ even or odd.
\end{proof}

While the essential surface just constructed has only one boundary component and hence is non-orientable, by Lemma~\ref{orient} there also exists an orientable candidate surface with twice as many sheets and boundary components associated to the same admissible edgepath system, so it too is essential with the same boundary slope but twice the Euler characteristic; that is:

\begin{cor}
For any integers $m,n>1$ with $n$ odd, the exterior of the Montesinos knot $K(-1/2, m/(2m\pm1),1/n)$ in $S^3$ contains an orientable essential surface with $2n$ sheets, two boundary components, Euler characteristic $-2n$, and boundary slope $2(n-1)^2/n$ if $m$ is even and $2(n+1)^2/n$ if $m$ is odd.
\end{cor}

We now prove Theorem \ref{res2} restated here for convenience.

\begin{thm}\label{reres2}
For any integer $m\geq1$ and odd integer $n \geq 4m + 1$, the exterior of the Montesinos knot $K=K(-1/2, m/(2m+1),1/n)$ in $S^3$ contains an essential surface with $n-4m+1$ sheets, two boundary components, Euler characteristic $2m+3-n$, and boundary slope $2(n^2-4mn-n+8m-1)/(n-4m+1)$ if $m$ is even and $2(n^2-4mn+3n-8m+3)/(n-4m+1)$ if $m$ is odd.
\end{thm}

\begin{proof}
Let $\gamma$ be the edgepath system given by
\begin{eqnarray*}
\gamma_1 &=& \left[\frac{n-4m+1}{n-2m}\sbrac{-\frac{1}{2}} + \frac{2m-1}{n-2m}\sbrac{-\frac{1}{2}}^\circ\right]\\
\gamma_2 &=& \left[\frac{n-4m-1}{n-4m+1}\sbrac{\frac{1}{2}} + \frac{2}{n-4m+1}\sbrac{\frac{m}{2m+1}}, \sbrac{\frac{m}{2m+1}}\right]\\
\gamma_3 &=& \left[\frac{n-4m}{n-4m+1}\sbrac{0} + \frac{1}{n-4m+1}\brac{1}{n}, \sbrac{\frac{1}{n}}\right]
\end{eqnarray*}

We verify that $\gamma$ satisfies conditions (E1-4):
\begin{itemize}
\item[(E1)] $\gamma_1$ is constant and lies on $\left[\sbrac{-1/2}, \sbrac{-1/2}^\circ\right]$, $\gamma_2(0) = \sbrac{m/(2m\pm1)}$ lies on $\left[\sbrac{m/(2m\pm1)}, \sbrac{m/(2m\pm1)}^\circ\right]$, and $\gamma_3(0) = \sbrac{1/n}$ lies on $\left[\sbrac{1/n}, \sbrac{1/n}^\circ\right]$.
\item[(E2)] No $\gamma_i$ stops and retraces itself or travels along two sides of a triangle in $\DD$ in succession.
\item[(E3)] Using Formula (\ref{coord}) and the formula that follows it,
\begin{eqnarray*}
\gamma_1(1) &=& \frac{n-4m+1}{n-2m}\sbrac{-\frac{1}{2}} + \frac{2m-1}{n-2m}\sbrac{-\frac{1}{2}}^\circ = \left(\frac{n-1}{2n-4m},-\frac{1}{2}\right)\\
\gamma_2(1) &=& \frac{n-4m-1}{n-4m+1}\sbrac{\frac{1}{2}} + \frac{2}{n-4m+1}\sbrac{\frac{m}{2m+1}} \\&=& \left(\frac{n-1}{2n-4m},\frac{n-2m-1}{2n-4m}\right)\\
\gamma_3(1) &=& \frac{n-4m}{n-4m+1}\sbrac{0} + \frac{1}{n-4m+1}\brac{1}{n} =\left(\frac{n-1}{2n-4m},\frac{1}{2n-4m}\right)
\end{eqnarray*}
These lie on a vertical line (i.e., have the same $u$-coordinate), and their $v$-coordinates sum to zero.
\item[(E4)] Each $\gamma_i$ proceeds monotonically from right to left.
\end{itemize}  

Hence, $\gamma$ is an admissible edgepath system, so a candidate surface $S$ can be associated to it. Because $\gamma_1$ is constant, $S$ is essential by Lemma~\ref{const}.

Now $|\gamma_2| = \frac{n-4m-1}{n-4m+1}$ and $|\gamma_3| = \frac{n-4m}{n-4m+1}$, so $n-4m+1$ is the least positive integer such that $n|\gamma_2|, n|\gamma_3| \in \Z$, and $(n-4m+1)\frac{n-2m}{n-4m+1} \in \Z$ for the constant edgepath $\gamma_1$, so $S$ has $n-4m+1$ sheets and Euler characteristic $2m+3-n$ by Formula (\ref{euler}).

Since $\gamma_1$ is constant, $\gamma_2$ is $\frac{n-4m-1}{n-4m+1}$ of an edge traveling upward, and $\gamma_3$ is $\frac{n-4m}{n-4m+1}$ of an edge traveling downward from right to left, $e_+=\frac{n-4m-1}{n-4m+1}$, $e_-=\frac{n-4m}{n-4m+1}$, and $\tau(S) = 2/(n-4m+1)$.

The same admissible edgepath systems whose associated candidate surfaces were Seifert surfaces $\Sigma$ for $K$ in the proof of Theorem~\ref{reres} work here as well; recall that $\tau(\Sigma) = 4-2n$ if $m$ is even and $-4-2n$ if $m$ is odd, so the boundary slope of $S$ is $\tau(S) - \tau(\Sigma) = 2(n^2-4mn-n+8m-1)/(n-4m+1)$ if $m$ is even and $2(n^2-4mn+3n-8m+3)/(n-4m+1)$ if $m$ is odd.

Because $n$ is odd, $n-4m+1$ is even, so the boundary slopes above are reducible by 2. Rewriting as $2n-4+\frac{2}{n-4m+1}$ if $m$ is even and $2n+4+\frac{2}{n-4m+1}$ if $m$ is odd shows that in either case the denominator of the boundary slope in lowest terms is $(n-4m+1)/2 \in \Z$, so since $S$ has $n-4m+1$ sheets, it has two boundary components.
\end{proof}

Since the essential surface just constructed has two boundary components, it may be orientable; e.g., if $m=1$ and $n=5$, $K(-1/2,1/3,1/5)$ is the $(-2,3,5)$-pretzel knot, and the essential surface is an annulus, and if $m=1$ and $n=7$, $K(-1/2,1/3,1/7)$ is the $(-2,3,7)$-pretzel knot, and the essential surface is a twice-punctured torus (cf.~\cite{Wu11}). Reducibility of the boundary slope was overlooked in \cite{Ichihara14}, which considers only the case $m=2$ (i.e., the knots $K(-1/2,2/5,1/n)$) and incorrectly concludes that the essential surface has only one boundary component and hence must be non-orientable; in this case, the denominator of the boundary slope is actually half the Euler characteristic of the essential surface, which in fact has two boundary components. In any case, orientable essential surfaces can always be obtained by Lemma~\ref{orient} as before.

Generalizing the main result of \cite{Ichihara14}, which treats only the case $m=2$, we see that infinitely many knots have boundary slopes of arbitrarily small difference:

\begin{cor}
For any $\epsilon > 0$ and integer $m>1$, there exists an integer $N$ such that for all $n \geq N$, the Montesinos knot $K(-1/2, m/(2m+1),1/n)$ has boundary slopes with difference less than $\epsilon$.
\end{cor}

\begin{proof}
Let $\epsilon > 0$ and $m > 1$. If $n \geq 4m + 1$, $K(-1/2, m/(2m+1),1/n)$ has boundary slopes $2(n-1)^2/n$ and $2(n^2-4mn-n+8m-1)/(n-4m+1)$ if $m$ is even and $2(n+1)^2/n$ and $2(n^2-4mn+3n-8m+3)/(n-4m+1)$ if $m$ is odd by Theorems~\ref{reres} and \ref{reres2}. In either case, their difference is $(8m-2)/(n^2-4mn+n)$, which decreases to 0 as $n$ increases. Thus, there exists an integer $N  \geq 4m + 1$ such that $(8m-2)/(n^2-4mn+n) < \epsilon$ for all $n \geq N$.
\end{proof}

Theorem~\ref{MSY} establishes that all knots considered in this section are (1,1)-knots, so Theorem~\ref{reres} shows that any odd integer $n>1$ occurs in the denominator of a boundary slope for infinitely many (1,1)-knots, and the proof of Theorem~\ref{reres2} shows that all nonzero integers occur in the denominator of a boundary slope for infinitely many (1,1)-knots:

\begin{cor}
For any integers $m, q \geq 1$, $K(-1/2, m/(2m+1),1/(4m+2q-1))$ has a boundary slope with denominator $q$.
\end{cor}

Note that throughout this paper we assume that denominators are positive, but boundary slopes are negated for a knot's mirror image.

\section{Discussion}\label{disc}
We conclude with discussion of the following bound on the denominators of boundary slopes for Montesinos knots (\cite[Theorem 1.1]{IchiharaMizushima07}).

\begin{thm}\label{upbd}
Let $K$ be a Montesinos knot with at least three rational tangles and $S$ an essential surface in the exterior of $K$ in $S^3$ with boundary slope $p/q$. If $K$ is not a $(-2,3,n)$-pretzel knot for odd $n \geq 3$, then $q \leq -\chi(S)/|\partial S|$. If $K$ is a $(-2,3,n)$-pretzel knot for odd $n \geq 3$, then $q \leq -\chi(S)/|\partial S| + 1$.
\end{thm}

By Theorem~\ref{reres}, for all integers $m,n>1$ with $n$ odd, the Montesinos knot $K(-1/2, m/(2m\pm1),1/n)$, which is not equivalent to any $(-2,3,t)$-pretzel knot, achieves the first bound. Only knots equivalent to the case $m=2$ are discussed in \cite{IchiharaMizushima07} and deemed ``the best possibility" (see \cite{Zieschang84} for equivalence of Montesinos knots). Thus, Theorem~\ref{reres} generalizes this result to all integers $m>1$.

If $m=1$ in Theorem~\ref{reres2}, we see that for any odd integer $n \geq 5$ the Montesinos knot $K(-1/2, 1/3,1/n)$, which is equivalent to the $(-2,3,n)$-pretzel knot, has boundary slope with denominator $(n-3)/2$ and hence achieves the second bound above. Knots equivalent to these are also discussed in \cite{IchiharaMizushima07} and deemed ``the best possibility."

If $m=2$ in Theorem~\ref{reres2}, we see that for any odd integer $n \geq 9$ the Montesinos knot $K(-1/2, 2/5, 1/n)$, which is not equivalent to any $(-2,3,t)$-pretzel knot, has boundary slope with denominator $(n-7)/2$ and hence achieves the first bound above. Knots equivalent to these are discussed in \cite{IchiharaMizushima07} but are not deemed ``the best possibility."

Finally, by Theorem~\ref{reres2}, for all $m>2$ and odd $n \geq 4m+1$, the Montesinos knot $K(-1/2, m/(2m+1), 1/n)$, which is not equivalent to any $(-2,3,t)$-pretzel knot, has boundary slope with denominator $(n-4m+1)/2$ and hence achieves neither bound. Knots equivalent to these are also mentioned in \cite{IchiharaMizushima07}.

\section*{Acknowledgments}
The author thanks Alan Reid for guidance that initiated this work. Nathan Dunfield's computer program to compute boundary slopes for Montesinos knots (available at {\tt www.CompuTop.org} and described in \cite{Dunfield01}) was used in formulating and checking cases of our results. The author also thanks Eric Chesebro for assistance with this program, the algorithm of \cite{HatcherOertel89}, and the figures.

\bigskip

\noindent\address{Department of Mathematics\\
  St.~Edward's University\\
  3001 South Congress Ave\\
  Austin, TX 78704, USA}
  
\smallskip  

\noindent\email{jasonc@stedwards.edu}

\end{document}